\documentclass[11pt]{article}
\usepackage{latexsym}
\usepackage{amsfonts}
\makeatletter

\def\@footnotetext#1{\insert\footins{%

\footnotesize

 \interlinepenalty\interfootnotelinepenalty

 \splittopskip\footnotesep

 \splitmaxdepth \dp\strutbox \floatingpenalty \@MM

 \hsize\columnwidth \@parboxrestore

 \edef\@currentlabel{\csname p@footnote\endcsname\@thefnmark}\@makefntext
 {\rule{\z@}{\footnotesep}\ignorespaces
#1\strut}}}

\def\abstract{\small\quotation{\hskip-\parindent\sc Abstract.}}
\def\classification{\@ifnextchar [{\@xfootnotenext}%
 {\begingroup\let\protect\noexpand
 \xdef\@thefnmark{}\endgroup
 \@footnotetext}}

\title {}

\begin{document}

\classification {{\it 2000 Mathematics Subject Classification:} Primary 14E09, 
 14E25; Secondary 14A10, 13B25.
$\ast$) Partially supported by RGC Grant Project 7126/98P.}

\begin{center}

{\bf \Large Birational morphisms of the plane}

\bigskip

{\bf Vladimir Shpilrain}

\medskip

 and

\medskip

 {\bf Jie-Tai Yu}$^{\ast}$

\end{center}

\medskip

\begin{abstract}
\noindent  Let $A^2$ be the affine plane over a field $K$ of characteristic $0$.
Birational morphisms of $A^2$ are mappings  $A^2 \to A^2$ given by polynomial 
 mappings $\varphi$ of the polynomial algebra $K[x,y]$ such that for the quotient 
fields, one has $K(\varphi(x), \varphi(y)) = K(x,y)$. 
 Polynomial automorphisms are obvious examples of such mappings.   
Another obvious example is the mapping $\tau_x$ given by $x \to x, ~y \to xy$. 
For a while, it was an open question whether every birational morphism is 
a product of polynomial automorphisms and copies of $\tau_x$. This question 
was answered in the negative by P. Russell (in an informal communication). 
 In this paper, we give a  simple combinatorial solution of the same problem. 
 More importantly, our method     yields 
an  algorithm for deciding 
whether a given birational morphism can be factored that way. 
\end{abstract}

\bigskip

\noindent {\bf 1. Introduction }
\bigskip

 Let $ K[x, y]$ be the polynomial algebra in two variables 
 over a field $K$ of characteristic $0$, and $A^2$  the affine plane over $K$.  
Birational morphisms of 
$A^2$ are mappings  $A^2 \to A^2$ given by polynomial 
 mappings $\varphi$ of the   algebra $K[x,y]$ such that for the quotient 
fields, one has $K(\varphi(x), \varphi(y)) = K(x,y)$. 
 Polynomial automorphisms are obvious examples of such mappings.   
Another obvious example is the mapping $\tau_x$ given by $x \to x, ~y \to xy$. 
Call either of those mappings a {\it simple affine contraction}. 
For a while, it was an open question whether every birational morphism is 
a product of simple affine contractions. This question 
was answered in the negative by  Russell (in an informal communication). 
The paper \cite{Daigle} by Daigle contains  an elaboration of Russell's 
methods.

 In this paper, we use  an  altogether different method of ``peak reduction" 
to easily establish the same result. More importantly, our method gives   
an  algorithm for deciding 
whether a given birational morphism can be factored that way:
\medskip

\noindent {\bf Theorem 1.1.} Let $K$ be an algebraically closed field.  
Then there is  an  algorithm for deciding 
whether a given birational morphism of $A^2$  over $K$ is a product of 
simple affine contractions.  
\medskip

Here we assume that we are able to perform calculations 
 in the ground field $K$, which basically means that, given 
(arithmetic expressions for) 
two elements of $K$, we can decide whether or not they are equal. 

 Our method can also be applied to fields that are not algebraically closed 
to prove that a particular birational morphism is {\it not} 
a product of simple affine contractions. 

\medskip

\noindent {\bf Proposition 1.2.} Let $K$ be any field of characteristic $0$,
 and $\phi : x \to  u,  ~y \to v$ 
a birational morphism  of $A^2$  over $K$. If $\phi$ is a product of simple affine 
contractions, then the sum of the  degrees of $u$ and $v$ can be decreased by a 
single generalized simple affine contraction. 
\medskip

 Generalized simple affine contractions are morphisms of the types (ET1)--(ET3), 
defined in Section 2. 

 For example, the following  birational 
morphism  of $A^2$  over ${\bf R}$  is not a product of simple affine 
contractions (see Example 2 in the end of Section 2): $x \to  x,  ~y \to  yx^2 + y$. 
 To find a birational morphism of $A^2$  over ${\bf C}$ that is not 
a product of simple affine contractions,  is somewhat more difficult. 
The following birational morphism is a   simplification of  an example 
due to    Cassou-Nogues and  Russell \cite{CN} : 

\noindent  $x \to x^4 y^2 -2x^3y+x^2+xy, \\
y \to  x^6y^3 -3x^5y^2 +3x^4 y + 2x^3y^2 - x^3 - 3x^2y +x +y$.

 We explain this example in Section 2, after the proof of Theorem 1.\\

\noindent {\bf 2. Peak reduction }
\bigskip

The ``peak reduction" method is a simple but rather powerful combinatorial 
technique with applications in many different areas of mathematics 
as well as theoretical computer science. It was introduced 
by Whitehead,  who used it
to solve an important algorithmic problem concerning automorphisms of a free
group (see e.g. \cite{Lyndonbook}). Since then, this method has been  used to solve 
various problems in
group theory,   topology, combinatorics, and probably in some other areas 
as well. 

In general, this method is used  to find  some kind of canonical 
form of a given object $P$ under the action of a given group (or a semigroup) 
$T$  of transformations. 
 The  idea behind the    method is rather simple: one chooses the 
{\it complexity} of  an object $P$ one way or another, and declares a
 canonical form of $P$ an object $P'$ whose complexity is minimal 
among all objects $t(P), ~t \in T$. To actually find a canonical form, or a 
 ``canonical model",  
 $P'$ of a given object $P$, one tries to arrange a sequence of sufficiently 
simple transformations so that the complexity of an object decreases 
{\it at every step}. To prove that such an arrangement is possible, one uses 
``peak reduction"; that means, if in some sequence of simple transformations 
the complexity goes up (or remains unchanged) before eventually going down, 
then there must be a pair of {\it consecutive} simple transformations in the 
sequence (a ``peak") such 
  that  one of them increases the  complexity  (or leaves it 
unchanged), and then the other one decreases it. Then one tries to prove 
that such a peak can always be reduced. 

 In the commutative algebra context, objects are polynomials; their complexity 
is their degree; the group of transformations is the group of polynomial 
automorphisms; simple transformations are elementary and 
linear  automorphisms. (An elementary automorphism is a one that changes 
just one variable.) 

 We have used this technique in our earlier paper \cite{ShYu} to 
contribute toward a classification 
 of two-variable polynomials having  classified, up to  an automorphism,  
polynomials of the form $ax^n + by^m + \sum_{im+jn \le
mn} c_{ij} x^i y^j$  (i.e.,  polynomials whose 
Newton polygon is  either a triangle or a line segment). Later, 
Wightwick \cite{Wightwick} used the idea of peak reduction in 
combination with  splice diagrams technique due to  Eisenbud and  Neumann  
\cite{EN} 
to classify {\it all} two-variable  polynomials over ${\bf C}$ up to  an
automorphism. More details and results can be found in our survey \cite{ShYusurvey}. 

 Here we use this technique to prove our  statements. 

\medskip 
 
\noindent {\bf Proof of Theorem 1.1 and Proposition 1.2.} 
 Consider the direct product $K[x,y] \times K[x,y]$ of two 
copies of the  polynomial algebra   $K[x,y]$, and introduce the following 
 elementary transformations (ET) that can be applied to elements of this 
 direct product: 
\smallskip 

\noindent {\bf (ET1)} $(u, ~v) \longrightarrow (\frac{u+a}{c\cdot v+b}, ~v)$ 
for arbitrary $a, b, c \in K$, with the denominator  not 0. 
\smallskip 

\noindent {\bf (ET1$'$)} $(u, ~v) \longrightarrow (u, ~\frac{v+a}{c\cdot u+b})$ 
for  arbitrary $a, b, c \in K$, with the denominator  not 0. 
\smallskip 

\noindent {\bf (ET2)} $(u, ~v) \longrightarrow (u + q(v), ~v)$ 
for an arbitrary polynomial $q(v)$. 
\smallskip 

\noindent {\bf (ET2$'$)} $(u, ~v) \longrightarrow (u, ~v+ q(u))$ 
for an arbitrary polynomial $q(u)$. 
\smallskip 

\noindent {\bf (ET3)} $(u, ~v) \longrightarrow (v, ~u)$. 

\smallskip

 Transformations (ET1) or (ET1$'$) are only applied if the corresponding 
ratio is a polynomial. 
\smallskip

 It is clear that   a birational morphism $x \to u, ~y \to v$
 is a product of polynomial 
automorphisms and  mappings  of the form  $x \to x, ~y \to x\cdot y$
 if and  only if the pair $(u, ~v)$ can be taken to  $(x, ~y)$ by a 
sequence of elementary transformations (ET1)--(ET3).

 We are now going to use ``peak reduction" to show that, if the 
sum of the  degrees  
of a pair  of polynomials can be reduced by a sequence of elementary 
transformations (ET1)--(ET3), then it can be reduced by a single elementary 
transformation. 

 Obviously, (ET3)  cannot change the sum of the  degrees; also, a  single 
 (ET1) or (ET1$'$) can only decrease the sum of the  degrees, unless $c=0$. 
 Thus, up to a symmetry, there are only 
the following possibilities for a ``peak":

\smallskip 

\noindent {\bf (1)} (ET3) followed by one of the (ET1), (ET1$'$), 
(ET2), or (ET2$'$). 

\smallskip 

\noindent {\bf (2)} (ET2) followed by (ET2$'$). 
\smallskip

\noindent {\bf (3)} (ET2) followed by (ET1) or (ET1$'$).  
\smallskip

 The first possibility is trivial: for example, (ET3) followed by (ET1) 
is the same as (ET1$'$) followed by (ET3). Therefore, if the 
sum of the  degrees  could be reduced by (ET3) followed by (ET1), then it 
could  also be reduced by (ET1$'$) followed by (ET3), hence by just (ET1$'$). 

 The second possibility was handled in \cite{ShYu}; it turns out to be 
easy, too. 

 In (3), if (ET2) is followed by (ET1$'$), we get the pair 
$(u+q(v),~\frac{v+a}{c\cdot [u+q(v)]+b})$. Since we are under the assumption 
that (ET2) does not decrease the sum of the  degrees, we have 
 $deg(u+q(v)) \ge deg(u)$. Therefore, we have two possibilities:
\smallskip 

\noindent {\bf (i)} $deg(u+q(v)) > deg(u)$. In this case, 
$deg(u+q(v)) = deg(q(v))= deg(q) \cdot  deg(v)$, and we conclude that 
 the polynomial $q$ 
must be linear; otherwise, $v+a$ cannot be divisible by $c\cdot [u+q(v)]+b$. 
Moreover, $q$ has to be a constant since otherwise, the ratio 
$~\frac{v+a}{c\cdot [u+q(v)]+b}$ is a constant, but no pair of the form  
$(u, ~c), ~c \in K$, can be taken to  $(x, ~y)$ by a sequence of  (ET1)--(ET3).
But if $q(v)$ is a constant, then (ET1$'$) cannot change the sum of the  degrees.
\smallskip 

\noindent {\bf (ii)} $deg(u+q(v)) = deg(u)$. In this case, 
$deg(u) \ge deg(q(v))$, so  again,   the polynomial $q$ 
must be linear; otherwise, $v+a$ cannot be divisible by $c\cdot [u+q(v)]+b$. 
As in (i) above, we conclude that $q$ has to be a constant, in which case 
 (ET1$'$) cannot change the sum of the  degrees.
\smallskip

 Thus, we are left with the crucial possibility where (ET2) is followed by (ET1)
 with $c \ne 0$. 
The result of applying this pair of elementary transformations to $(u, v)$ 
is $(\frac{u+q(v)+a}{c\cdot v+b}, ~v)$. The polynomial $q(v)$ can be written as 
$r(c\cdot v+b)$ for some other polynomial $r$. Then 
$$(\frac{u+q(v)+a}{c\cdot v+b}, ~v) =
 (\frac{u+r(c\cdot v+b)+a}{c\cdot v+b}, ~v) = 
(\frac{u+r_1(c\cdot v+b)+a_1}{c\cdot v+b}, ~v),$$

\noindent  where the polynomial  $r_1$ has zero constant term. 

 Since $r_1(c\cdot v+b)$  is divisible by $c\cdot v+b$, we get 
$$(\frac{u+r_1(c\cdot v+b)+a_1}{c\cdot v+b}, ~v) = 
(\frac{u+a_1}{c\cdot v+b}+r_2(c\cdot v+b), ~v).$$

 This means, in particular, that $u+a_1$  is divisible by $c\cdot v+b$. But then 
a single (ET1) would reduce the sum of the  degrees of the pair $(u, v)$. 
This completes the ``peak reduction".

 To get now an  algorithm claimed in the statement of Theorem 1, 
  we first  have to    show that one can effectively   determine 
whether a single elementary transformation can reduce the 
sum of the  degrees  of a given pair  of polynomials. For (ET2) or (ET2$'$), 
this  is well known (see e.g. \cite[Theorem 6.8.5]{PMCohn}). 
For (ET1) or (ET1$'$), the procedure is quite straightforward. Suppose 
we want to find out whether, say, $u(x,y)+a$  is divisible by $v(x,y)+b$ 
for some $a, b \in K$ (we can clearly assume that $c=1$).  We apply the usual 
``long division" algorithm 
with respect to some fixed term ordering (see e.g. 
\cite{AL}). In the end, the divisibility 
condition translates into a system of polynomial equations over $K$ in a 
single  variable $b$, complemented by an equation  $a=p(b)$ for some 
 polynomial $p$.  The solvability of such a system over an algebraically 
closed field can be decided by using Gr\"{o}bner bases technique (see e.g. 
\cite{AL}). (Of course, an actual solution may not be found in general.) 

 If the system has no solutions, then a single elementary transformation cannot
 reduce the 
sum of the  degrees  of a given pair  of polynomials, and we conclude that 
our birational morphism  is not a product of 
simple affine contractions. If the system has a solution, then we apply 
the corresponding elementary transformation, keeping our parameters $a, b$, 
etc. as variables, because we do not have explicit values for them. 
 Then, for the new pair of polynomials, we do the same thing: we check if 
a single elementary transformation can reduce the 
sum of the  degrees. This will yield another system of polynomial equations over $K$,
that expands the previous system by   introducing new variables and new 
equations. Again, it can be decided whether or not this system has a solution. 

 This procedure  will  obviously terminate in a number of steps not exceeding 
the sum of the  degrees  of the   original pair  of polynomials.
 This completes the proof. $\Box$  
\medskip 
 
\noindent {\bf Example 1.} The birational morphism 

\noindent  $x \to u(x,y) = x^4 y^2 -2x^3y+x^2+xy, \\
y \to  v(x,y) = x^6y^3 -3x^5y^2 +3x^4 y + 2x^3y^2 - x^3 - 3x^2y +x +y$, 

\noindent mentioned in the  Introduction, is constructed as follows.
Start with the birational morphism  $(xy+1, ~x^2y+x)$, ~apply $x \to x+y, 
~y \to y$ to get $((x+y)y+1, ~(x+y)^2y+(x+y))$. Then apply 
$x \to \frac{1}{x}, ~y \to -x+x^2y$  (this is a bijective  morphism of $K(x,y)$) 
to get the above example. 

 A single transformation (ET2$'$) cannot reduce the 
sum of the  degrees  of the   pair $(u, ~v)$  because the  degree of 
$v$ is not divisible by the  degree of $u$ (cf. \cite{ShYu}). 
 Neither can a single (ET1$'$) : a straightforward check shows that 
$v(x,y)+a$  is not  divisible by $c\cdot u(x,y)+b$ for any $a, b, c \in K$.
Therefore, by Proposition 1.2, this birational morphism 
is not a product of simple affine contractions. 
\medskip 
 
\noindent {\bf Example 2.} The birational morphism $x \to  x,  ~y \to  yx^2 + y$ 
is not a product of simple affine contractions over ${\bf R}$, by Proposition 1.2. 
However, over ${\bf C}$, a single (ET1$'$)  can reduce the degree of $yx^2 + y$ 
since $yx^2 + y = y(x^2 + 1)=y(x+i)(x-i)$. It is now easy to see that  over ${\bf C}$,
 this birational morphism is a product of simple affine contractions. \\

\noindent {\bf Acknowledgement}
\medskip

\indent We are grateful to Andrew Campbell   for  helpful comments.

\noindent 
 Department of Mathematics, The City  College  of New York, New York, 
NY 10031 
\smallskip

\noindent {\it e-mail address\/}: ~shpil@groups.sci.ccny.cuny.edu 
\smallskip

\noindent {\it http://www.sci.ccny.cuny.edu/\~\/shpil/} \\

\noindent Department of Mathematics, The University of Hong Kong, 
Pokfulam Road, Hong Kong 
\smallskip

\noindent {\it e-mail address\/}: ~yujt@hkusua.hku.hk 
\smallskip

\noindent {\it http://hkumath.hku.hk/\~\/jtyu}


\baselineskip 11 pt


 \begin{thebibliography}{ABC}



\bibitem{AL}
W.~Adams and P.~Loustaunau,
{\sl An introduction to Gr\"{o}bner bases\/}.
American Mathematical Society, Providence, 1994. 


\bibitem{CN}
P. Cassou-Nogues and  P. Russell, {\it On some birational 
endomorphisms of the affine plane}, preprint.  


\bibitem{PMCohn}
  P.~M.~Cohn, {\sl Free Rings and Their Relations. 2nd Ed.\/} 
  Academic Press, London, 1985.


\bibitem{Daigle}
D. Daigle,  {\it Birational endomorphisms of the affine plane}, 
J. Math. Kyoto Univ. {\bf 31} (1991),   329--358.

\bibitem{EN} 
D. Eisenbud and W. D. Neumann,
  {\it   Three-dimensional link theory and invariants of plane curve
    singularities.}  Ann. Math. Stud. {\bf 110}, Princeton.  Princeton
  Univ. Press (1985).


\bibitem{Lyndonbook}
R. Lyndon, P. Schupp, {\sl Combinatorial Group Theory}, 
(Reprint of the 1977 edition). In {\sl Classics in Mathematics},
Springer-Verlag, Berlin (2001).


\bibitem{Wightwick}
P. Wightwick, {\it  Equivalence of polynomials 
under automorphisms}, J. Pure Appl. Algebra {\bf 157} (2001), 341--367. 

\bibitem{ShYu}
V. Shpilrain and J.-T. Yu,  {\it Embeddings 
of curves in the plane},  J. Algebra {\bf 217} (1999),  668--678. 


\bibitem{ShYusurvey}
V. Shpilrain and J.-T. Yu,  {\it Peak reduction 
technique in commutative algebra: a survey}, in:  Combinatorial and computational 
algebra (Hong Kong, 1999), 237--247, Contemp. Math. {\bf  264}, Amer.
Math. Soc., Providence, RI, 2000.\\












\end{thebibliography}
\end{document}